\documentclass{amsart}
\usepackage{amsaddr}
\usepackage[english]{babel}
\usepackage{latexsym}
\usepackage{amsmath,amssymb,amsfonts,amsthm}
\usepackage{mathrsfs}

\title{Foliations on $\mathbb{CP}^3$ of degree $2$ that have a line as singular set}
\author[]{Claudia R. Alc\'antara$^1$, Dominique Cerveau$^2$}

\address{$^1$Departamento de Matem\'aticas, Universidad de Guanajuato.
 Callej\'on Jalisco s/n, C.P. 36000,
Guanajuato, Gto. M\'exico.}
\email{claudia@cimat.mx}
\address{$^{1,2}$Universit\'e de Rennes I, IRMAR, Campus Beaulieu, 35042, Rennes Cedex, France.}
\email{dominique.cerveau@univ-rennes1.fr}

 \thanks{This work was supported by CNRS, UMR6625, INSTITUT DE RECHERCHE
MATHEMATIQUE DE RENNES (IRMAR)}

\newcommand{\C}{\mathbb{C}}

\newcommand{\Q}{\mathbb{Q}}
\newcommand{\Z}{\mathbb{Z}}
\newcommand{\N}{\mathbb{N}}

\newcommand{\CP}{\mathbb{CP}^2}

\newtheorem{teo}{Theorem}
\newtheorem{defin}{Definition}

\newtheorem{cor}{Corollary}
\newtheorem{ex}{Example}

\begin{document}
\maketitle

\begin{abstract}
In this work we classify foliations on $\mathbb{CP}^3$  of codimension 1 and degree $2$ that have a line as singular set. To achieve this, we do a complete description of the components. We prove that the boundary of the exceptional component has only 3 foliations up to change of coordinates, and this boundary is contained in a logarithmic component. Finally we construct  examples of foliations on $\mathbb{CP}^3$  of codimension 1 and degree $s \geq 3$ that have a line as singular set and such that 
they form a family with a rational first integral of degree $s+1$ or they are logarithmic foliations where some of them have a minimal rational first integral of degree not bounded.
\end{abstract}

 \section{Introduction}

 A holomorphic foliation $\mathcal{F}$ on the projective space $\mathbb{CP}^3$ of codimension one and degree $s$  is given by the projective class of a 1-form:
 
$$ \omega=A_1(z_1,z_2,z_3,z_4)dz_1+A_2(z_1,z_2,z_3,z_4)dz_2+A_3(z_1,z_2,z_3,z_4)dz_3+A_4(z_1,z_2,z_3,z_4)dz_4,$$

\noindent where $A_1, A_2, A_3, A_4 \in \C[z_1,z_2,z_3,z_4]$ are homogeneous of degree $s+1$, and they satisfy:

\begin{enumerate}
\item $\sum_{i=1}^4 z_iA_i(z_1,z_2,z_3,z_4)=0$
\item The integrability condition:  $\omega \wedge d\omega=0$.
\end{enumerate}

\noindent We consider the subspace of classes of 1-forms that satisfy these conditions and such that its singular set $Sing(\omega)$ has codimension 2, we denote this space by $\mathcal{F}(s,3)$. Then $\mathcal{F}(s,3)$ can be identified with a Zariski's open set in the projective space $\mathbb{P}H^0(\mathbb{CP}^3, \Omega_{\mathbb{CP}^3}^1(s+2))$, which has dimension $4 \binom{s+4}{3} - \binom{s+5}{3}-1$. During this article, when we talk about closure we will be thinking with respect to $\mathcal{F}(s,3)$, unless otherwise specified.
\\

\noindent The form $ \omega=A_1(z_1,z_2,z_3,z_4)dz_1+A_2(z_1,z_2,z_3,z_4)dz_2+A_3(z_1,z_2,z_3,z_4)dz_3+A_4(z_1,z_2,z_3,z_4)dz_4$ will be called a homogeneous expression of the foliation $\mathcal{F}$. The singular set of the foliation is the algebraic variety

$$Sing(\mathcal{F})=\mathbb{V}(A_1,A_2,A_3,A_4) \subset \mathbb{CP}^3,$$

\noindent which has dimension 1. Then the simplest singular set that a foliation on $\mathbb{CP}^3$ can have is a line, which in this case is the intersection of two hyperplanes. 
Due to the nature of this singular set it is clear that we can perform more calculations and try to prove more conjectures in the subspace of these foliations. It is for this reason that it is interesting to classify this type of foliations. 
Here we note that the classification of these foliations is also the question proposed on page 57 of the book  \cite{Deserti-Cerveau}.
\\

Then the main objective of this work is to find all the foliations in $\mathcal{F}(2,3)$ that have a line as singular set. Of course the linear pull-back 
of foliation on $\CP$ of degree $2$ with a unique singular point satisfy this property. We can see in the article \cite{Cerveau} that there exist just four foliations, up to change of coordinates, with one singular point, they are:

\begin{align*}
\nu_1&=z_3^3dz_1-z_1z_2^2dz_2+(z_2^3-z_1z_3^2)dz_2\\
\nu_2&=-z_3^3dz_1+z_3(z_2^2+z_1z_3)dz_2+(z_1z_3^2-z_2^3-z_1z_2z_3)dz_3\\
\nu_3&=z_2^2z_3dz_1-z_3(z_2^2+z_3^2)dz_2+(z_2^3+z_2^2z_3-z_1z_2z_3)dz_3\\
\nu_4&=(z_1z_2z_3+z_3^3-z_2^3)dz_1+z_2(z_3^2-z_1z_2)dz_2-(z_2^2z_3+z_1^2z_2+z_1z_3^2)dz_3.
\end{align*}

The first one has a rational first integral of degree $3$, and for the last one the singularity is a saddle-node and the foliation does not have algebraic invariant curves. With this in mind we can ask:
It will be that all the foliations on $\mathbb{CP}^3$ of degree $2$ that have a line as singular set are linear pull-back of these $4$ foliation on $\CP$? The answer is not, however all of them are somehow obtained from these $4$.
\\

More specifically, in this set there are, up to change of coordinates: the linear pull-back  of the three foliations $\nu_2$, $\nu_3$, $\nu_4$ on $\CP$ of degree $2$ with a unique singular point and without rational first integral and there are three foliations on 
$\mathbb{CP}^3$ with a rational first integral of degree $3$ where one of them is the linear pull-back of the foliation $\nu_1$ on $\CP$ (see Theorem \ref{Classification}).
 \\
 
 For the proof of this result we study each irreducible component of the space $\mathcal{F}(2,3)$ of foliations on $\mathbb{CP}^3$  of degree $2$. We know by the article  \cite{Cerveau-LinsNeto}  that $\mathcal{F}(2,3)$  has $6$ irreducible components,
 they are:

\begin{enumerate}
\item $S(2,3)$: the foliations which are linear pull-back of foliations on $\CP$ of degree $2$ with isolated singularities.
\item $\overline{R(2,2)}$: the closure of the space of foliations with a rational first integral $\frac{f(z_1,z_2,z_3,z_4)}{g(z_1,z_2,z_3,z_4)}$, where $f$ and $g$ in $\C[z_1,z_2,z_3,z_4]$ have degree 2 and $f$ defines a smooth hypersurface.
\item $\overline{R(1,3)}$: the closure of the space of foliations with a rational first integral $\frac{f(z_1,z_2,z_3,z_4)}{L^3(z_1,z_2,z_3,z_4)}$, where $f$ has degree $3$ and $L$ has degree $1$.
\item $\overline{L(1,1,1,1)}$: the closure of the space of logarithmic foliations given by the 1-forms:

$$\omega=\sum_{i=1}^4 \lambda_i L_1L_2L_3L_4 \frac{dL_i}{L_i},$$

\noindent where $L_1, L_2, L_3$ and $L_4$ in $\C[z_1,z_2,z_3,z_4]$ have degree $1$ and $\sum_{i=1}^4 \lambda_i=0$.

\item $\overline{L(1,1,2)}$: the closure of the space of logarithmic foliations given by the 1-forms:

$$\omega=\sum_{i=1}^3 \lambda_i f_1f_2f_3 \frac{df_i}{f_i},$$

\noindent where $f_1, f_2$ define different hyperplanes, $f_3$ defines an irreducible hypersurface of degree $2$ and $\lambda_1+\lambda_2+2\lambda_3=0$.

\item $\overline{E(3)}$: the exceptional component, which is the closure of the orbit of one foliation which singular set contains the twisted cubic.

\end{enumerate}

The components (4) and (5) are called logarithmic components.
\\

The most detailed and difficult analysis is that referring to the boundaries of the components.  In this sense we have obtained an interesting result that says that in the boundary of the exceptional component $\overline{E(3)}$ there exist, up to change of coordinates, only three foliations which are also in the logarithmic component $\overline{L(1,1,2)}$ (There are also three forms whose singular set has dimension 2, but for the sake of completeness we will also describe them). The proof of this result uses Geometric Invariant Theory and the techniques developed by Kirwan in the book \cite{Kirwan}. For the study of the boundary of the logarithmic components we mainly use the results of  \cite{Cerveau-Mattei}.
\\

The structure of the paper is as follows. In section 2 we describe all the foliations in the exceptional component, up to linear change of coordinates. Section three is devoted to study the logarithmic components. Finally in the last section we have the theorem of classification of foliations in $\mathcal{F}(2,3)$ that have a line as singular set. We finish with the construction of two examples of foliations in $\mathcal{F}(s,3)$, for $s \geq 3$ with a line as singular set and with one of the following properties:

\begin{enumerate}
\item They form a family that have a rational first integral of degree $s+1$.
\item They form a family of logarithmic foliations, where some of them have a minimal rational first integral whose degree is not bounded. 
\end{enumerate}

\section{The boundary of the exceptional component $\overline{E(3)}$} 

 We know that the exceptional component $\overline{E(3)}$ is the closure of the orbit of a foliation on $\mathbb{CP}^3$ which singular set contains the twisted cubic (see example 6 of \cite{Cerveau-LinsNeto}), this orbit is considering with respect to the action by the automorphism group 
 of $\mathbb{CP}^3$. Then to find the foliations that have a line as singular set it is necessary to find the foliations in the boundary of this orbit.
 For that we use Geometric Invariant Theory applied to the action by change of coordinates in this algebraic variety. We obtain the following.

 \begin{teo} \label{excepcional} Let $\overline{E(3)}$ be the exceptional component of the space of foliations on $\mathbb{CP}^3$ of codimension $1$ of degree $2$ then, up to change of coordinates, it contains only the foliations associated to the following 1-forms:
 
 \begin{enumerate}
 \item $\omega=(2z_2^2z_4-z_2z_3^2-z_1z_3z_4)dz_1+(2z_3^2z_1-3z_1z_2z_4)dz_2+(3z_1^2z_4-z_1z_2z_3)dz_3+(z_2^2z_1-2z_3z_1^2)dz_4$
\item $\omega_1+\omega_2$
\item $\omega_2+\omega_3$
\item $\omega_1+\omega_3$.
 \end{enumerate}

\noindent where 

\begin{align*}
&\omega_1=z_1(-z_3z_4dz_1+3z_1z_4dz_3-2z_3z_1dz_4),\\
&\omega_2=z_2(2z_2z_4dz_1-3z_1z_4dz_2+z_2z_1dz_4)\\
&\omega_3=z_3(-z_2z_3dz_1+2z_3z_1dz_2-z_1z_2dz_3).
\end{align*}
 \end{teo}
 
 \begin{proof} In the article \cite{Cerveau-LinsNeto} we can see that the exceptional component $\overline{E(3)}$ of the space of foliations on  $\mathbb{CP}^3$ of codimension $1$ of degree $2$  is the closure $\overline{SL_4(\C)\cdot \omega}$ in 
 $\mathcal{F}(2,3)$ of the orbit of the 
 foliation given by:

 \begin{align*}
 \omega&=(2z_2^2z_4-z_2z_3^2-z_1z_3z_4)dz_1+(2z_3^2z_1-3z_1z_2z_4)dz_2\\
&+ (3z_1^2z_4-z_1z_2z_3)dz_3+(z_2^2z_1-2z_3z_1^2)dz_4,
 \end{align*}
 
 \noindent with respect to the action by change of coordinates by the reductive algebraic group $SL_4(\C)$ on $\overline{E(3)}$.
 \\
 
 We can do this because we know that $\overline{E(3)}$ is an algebraic irreducible, variety of dimension $13$, therefore the linear action of $SL_4(\C)$  in the projective space $\mathbb{P}H^0(\mathbb{CP}^3, \Omega_{\mathbb{CP}^3}^1(4))=\mathbb{CP}^{44}$,  induces an action in this variety, which is 
 of course $SL_4(\C)$-invariant:
  
 \begin{align*}
 SL_4(\C) \times \overline{E(3)} &\to \overline{E(3)}\\
 (g,\nu) &\mapsto g \cdot \nu.
 \end{align*}
 
 \noindent We recall that a 1-parameter subgroup of an algebraic group $G$ is an algebraic morphism from $\C^*$ to $G$. And a very known result says that every 1-parameter subgroup of $SL_4(\C)$ is diagonalizable. Let:
 
  \begin{align*}
\lambda_{(n_1,n_2,n_3)}: \C^* &\to SL_4(\C)\\
 t  &\mapsto \left(\begin{array}{cccc}
t^{n_1}&0&0&0 \\
0&t^{n_2}&0&0\\
0&0&t^{n_3}&0\\
0&0&0&t^{n_4}
\end{array}\right),
\end{align*}

\noindent be a diagonal 1-parameter subgroup of $SL_4(\C)$, where $n_1, n_2, n_3, n_ 4 \in \Z$ and $n_1+n_2+n_3+n_4=0$. The action of $\lambda_{(n_1,n_2,n_3)}$ on $\omega$ is:

\begin{align*}
\lambda_{(n_1,n_2,n_3)}(t) \cdot \omega&=(2t^{n_3-n_4}z_2^2z_4-t^{n_2-n_3}z_2z_3^2-t^{n_1-n_2}z_1z_3z_4)dz_1\\
&+(2t^{n_3-n_4}z_3^2z_1-3t^{n_2-n_3}z_1z_2z_4)dz_2\\
&+(3t^{n_1-n_2}z_1^2z_4-t^{n_3-n_4}z_1z_2z_3)dz_3\\
&+(t^{n_2-n_3}z_2^2z_1-2t^{n_1-n_2}z_3z_1^2)dz_4.
\end{align*}

Let:
 
\begin{align*}
\omega_1&:=-z_1z_3z_4dz_1+3z_1^2z_4dz_3-2z_3z_1^2dz_4=z_1(-z_3z_4dz_1+3z_1z_4dz_3-2z_3z_1dz_4)\\
\omega_2&:=2z_2^2z_4dz_1-3z_1z_2z_4dz_2+z_2^2z_1dz_4=z_2(2z_2z_4dz_1-3z_1z_4dz_2+z_2z_1dz_4)\\
\omega_3&:=-z_2z_3^2dz_1+2z_3^2z_1dz_2-z_1z_2z_3dz_3=z_3(-z_2z_3dz_1+2z_3z_1dz_2-z_1z_2dz_3),
\end{align*}

\noindent then $\omega=\omega_1+\omega_2+\omega_3$ and

\begin{align*}
\lambda_{(n_1,n_2,n_3)}(t) \cdot \omega_1&=t^{n_1-n_2}\omega_1\\
\lambda_{(n_1,n_2,n_3)}(t) \cdot \omega_2&=t^{n_2-n_3}\omega_2\\
\lambda_{(n_1,n_2,n_3)}(t) \cdot \omega_3&=t^{n_3-n_4}\omega_3.\\
\end{align*}

\noindent If we take

\begin{align*}
\lambda_1(t)&=\lambda_{(3,-1,-1)}(t),\\ 
\lambda_2(t)&=\lambda_{(1,1,-1)}(t),\\ 
\lambda_3(t)&=\lambda_{(1,1,1)}(t), 
\end{align*}

\noindent then we have

\begin{equation*}
\lim_{t \to \infty} \lambda_i(t) \cdot \omega=\omega_i \quad \textrm{and} \quad  \lim_{t \to 0} \lambda_i(t) \cdot \omega=\omega_j+\omega_k
\end{equation*}

\noindent where $i,j,k \in \{1,2,3\}$ are $i,j,k$ different from each other.
\\

Then the classes of $\omega_1+\omega_2$, $\omega_2+\omega_3$ and $\omega_2+\omega_3$ are in 

$$\overline{SL_4(\C)\cdot \omega}-SL_4(\C) \cdot \omega=\overline{E(3)}-SL_4(\C)\cdot \omega,$$ 

\noindent here it is important to note that these are integrable 1-forms.
We note also that the 1-forms $\omega_1$, $\omega_2$, $\omega_3$ have a hyperplane as singular set and they are in the closure of $SL_4(\C)\cdot \omega$ but considered in $\mathbb{P}H^0(\mathbb{CP}^3, \Omega_{\mathbb{CP}^3}^1(4))$.
\\

In an analogous way we can see that for $i,j=1,2,3$, $i \neq j$, the classes of $\omega_i $ are in  $\overline{SL_4(\C) \cdot (\omega_i+\omega_j)}^{\mathbb{P}H^0(\mathbb{CP}^3, \Omega_{\mathbb{CP}^3}^1(4))}$ and the orbit $SL_4(\C)\cdot \omega_i$ is closed in $\mathbb{P}H^0(\mathbb{CP}^3, \Omega_{\mathbb{CP}^3}^1(4))$ because $\omega_i$ 
is the pull-back of a 1-form of degree $1$ of $\CP$ with three proper values and a line of singularities. Here we recall that the orbit of a matrix by conjugation is closed if the Jordan blocks have the minimum size, which is the case for three different proper values.
\\

\noindent On the other hand we have that $\lambda_{(3,1,-1)}(t) \cdot \omega=t^2\omega$, then 

\begin{equation*}
\lim_{t \to 0} \lambda_{(3,1,-1)}(t) \cdot \overline{\omega}=\lim_{t \to 0} t^2\overline{\omega}=0,
\end{equation*}

\noindent  where $\overline{\omega}$ is a different point from $0$ on $\omega$ in the affine cone of $\overline{E(3)}$. Then by Hilbert-Mumford criterion of 1-parameter subgroups (see \cite{Mumford}) the point $\omega$
is unstable for the action. We conclude that the set of semistable points is empty and the variety $\overline{E(3)}$ is the closed set of unstable points. Moreover, we have that $\{\lambda_{(3,1,-1)}(t): t \in \C^*\} \subset Aut(\omega)$, in fact we know that $\dim Aut(\omega)=2$. 
\\

Now we are going to prove that the above are all the foliations in $\overline{E(3)}$. For this we use theorem 12.26 that we can find in part II of the book  \cite{Kirwan}.
In order to simplify the arguments we are going to use the same notation as Kirwan in the mentioned book.
\\
 
 The theorem says that there exists a stratification by subvarieties $S_0,...,S_n$, of $\overline{E(3)}$ locally closed, invariant by the action, disjoint  and these subvarieties are parametrized by a finite set $\mathcal{B}=\{\beta_0, \beta_1,...,\beta_n\}$ of virtual 1-parameter subgroups of 
 $SL_4(\C)$ and:
 
 \begin{align*}
 S_0&=\overline{E(3)}^{ss}=\emptyset, \quad \textrm{semistable points for the action}\\ 
 \overline{E(3)}&=\overline{E(3)}^{un}=\bigcup_{i=1}^n S_{i} \quad \textrm{unstable points for the action}.
 \end{align*}
 
 \noindent Without giving too many details we are going to describe the construction of the stratum $S_i$. For that we remember that a diagonal 1-parameter subgroup of $SL_4(\C)$ is identified with a point in $\Z^3$ and that the set of virtual 1-parameter subgroups of $SL_4(\C)$ is identified with $\Q^3$.
\\

\noindent As we can see in definition $12.8$ of \cite{Kirwan} the indexing set $\mathcal{B}$ is the set of minimal combinations of weights lying in some Weyl chamber of the representation. And for $\beta_i \in \mathcal{B}$, we define:
  
 \begin{align*}
 Z_i&=\{\nu \in \overline{E(3)}: \beta_i(t) \cdot \nu=\nu\}\\
 Y_i&=\{\mu \in \overline{E(3)}: p_i(\mu)=\lim_{t \to 0} \beta_i(t) \cdot \mu \in Z_i\}
 \end{align*}

The function $p_i:Y_i \to Z_i$ is a locally trivial fibration with affine fiber  and  $p_i(\mu) \in \overline{SL_4(\C) \cdot \mu} \cap Z_i$ for all $\mu \in Y_i$. 
Finally $Y^{ss}_i=p_i^{-1}(Z_i^{ss})$, where $Z_i^{ss}$ is the set of semistable points with respect to a certain action. Then the stratum is $S_i=SL_4(\C) \cdot Y_i^{ss}$.
\vspace{.5cm}

On the other hand, it is easy to prove that the unique diagonal one parameter subgroup, up to integer multiples, which leaves $\omega$ fixed is $\beta_1=\lambda_{(3,1,-1)}$. Then with $\beta_1$ we will construct the stratum $S_1$, and it satisfies:
 
 $$SL_4(\C) \cdot \omega \subset S_1=SL_4(\C) \cdot Y^{ss}_1 \subset \overline{S_1} \subset \overline{E(3)},$$
 
 \noindent  then $\overline{S_1} = \overline{E(3)}$. In this case we have that the closed sets $Z_2,...,Z_n$ have intersection with  $Z_1$, because if this doesn't happen then we would have two foliations $\nu_1$, $\nu_2 \in \overline{E(3)}$  such that 
 $\overline{SL_4(\C) \cdot \nu_1} \cap \overline{SL_4(\C) \cdot \nu_2} \neq \emptyset$. Therefore we conclude that $\overline{S_1} = SL_4(\C) \cdot Y_1$. This means that it is enough to study the subvariety $Y_1$, but we will see that in fact it is enough to find the foliations in $Z_1$ to have all the foliations in $\overline{E(3)}$.
 \\
 
   The stabilizer of $\beta_1$ with respect to the adjoint action of $SL_4(\C)$ is the subgroup $D$ of diagonal matrices in $SL_4(\C)$. And the parabolic subgroup $P_1$ associated to $\beta_1$ is the group of upper triangular matrices (see definition 12.11 of \cite{Kirwan}).  Following page 153 of \cite{Kirwan}, we have that $D$ acts in $Z_1:=\{\nu \in \overline{E(3)}: \beta_1(t) \cdot \nu=\nu \}$ and $P_1$ acts in $Y_1$ therefore $D\cdot \omega$ and $P_1 \cdot \omega$ are open and dense in $Z_1$ and in $Y_1$, respectively.
\\ 
 
 Since $Z_1$ is irreducible we have that $Z_1=\overline{D \cdot \omega}$, and it will be enough to find the orbits in this closed set. Let $a_1,a_2,a_3,a_4 \in \C^{*}$ such that $a_1a_2a_3a_4=1$, then 
\begin{align*}
 \left(\begin{array}{cccc}
a_1&0&0&0 \\
0&a_2&0&0\\
0&0&a_3&0\\
0&0&0&a_4
\end{array}\right) \cdot \omega=a_1a_2^{-1}\omega_1+a_2a_3^{-1}\omega_2+a_3a_4^{-1}\omega_3.
\end{align*}

\noindent Thus $D \cdot \omega=\{\alpha_1 \omega_1+ \alpha_2 \omega_2+ \alpha_3 \omega_3: \alpha_1, \alpha_2, \alpha_3 \in \C^*\} $, note also  that $\omega_1$, $\omega_2$, $\omega_3$ 
are equivalence classes in the projectivization of different eigenspaces for the associated representation given by the action, with respect to the maximal torus $D$. We conclude that $\{\alpha_1 \omega+ \alpha_2 \omega_2+ \alpha_3 \omega_3: \alpha_1, \alpha_2, \alpha_3 \in \C\} \cap \overline{E(3)}$ is closed because is the intersection of the variety with a projective space, since this is contained in $Z_1$ it must be $Z_1$ and therefore:

$$Z_1=D\cdot\omega \cup \bigcup_{k\neq j}D\cdot (\omega_k+\omega_j).$$
\\

In this case the indexing set of virtual 1-parameter subgroups for constructing the stratification is 

$$\mathcal{B}=\{\beta_1=(3,1,-1), \beta_{12}=\lambda_3=(1,1,1), \beta_{13}=\lambda_2=(1,1,-1), \beta_{23}=\lambda_1=(3,-1,-1)\}.$$ 

\noindent Since $Y_1=\overline{P_1 \cdot \omega}$ we have the following sets and locally trivial fibrations:

\begin{align*}
&Y_1^{ss}=P_1 \cdot \omega,   &Y_{jk}^{ss}=P_{jk} \cdot (\omega_j+\omega_k)\\
& \downarrow p_1      \quad \quad  \quad  \quad \quad & \downarrow p_{jk} \quad \quad \quad \quad \quad \quad \\
&Z_1^{ss}=D \cdot \omega,  &Z_{jk}^{ss}=D_{jk}\cdot (\omega_j+\omega_k)
\end{align*}

\noindent  where $P_{jk}$ is the parabolic subgroups associated to $\beta_{jk}$ and $D_{jk}$ is the stabilizer of $\beta_{jk}$ with respect to the adjoint action (see page 154 of \cite{Kirwan}). Since $\overline{E(3)}=\overline{S_1}=SL_4(\C) Y_1$  we have that the unique orbits in the exceptional component, up to change of coordinates, are $\omega$, $\omega_1+\omega_2$, $\omega_2+\omega_3$ and $\omega_1+\omega_3$.
 
 \end{proof}
 
 Now we are going to see the singular set and the rational first integrals of the foliations in the exceptional component $\overline{E(3)}$. We use the notation of the above theorem.
 The singular set of $\omega$ is the union of the following three curves in $\mathbb{CP}^3$  
 \begin{enumerate}
 \item The conic in a plane $\mathbb{V}(z_1,2z_2z_4-z_3^2)$
 \item  The line $\mathbb{V}(z_1,z_2)$
 \item The twisted cubic $\mathbb{V}(2z_3^2-3z_2z_4,3z_1z_4-z_2z_3,z_2^2-2z_1z_3),$
 \end{enumerate}
 
 \noindent and its rational first integral is $\frac{(3z_4z_1^2-3z_1z_2z_3+z_2^3)^2}{(2z_1z_3-z_2^2)^3}$. For the other foliations we have: 
 
 \begin{align*}
 \textrm{Foliation} \quad& \textrm{Singular Set} & \textrm{Rational First Integral}\\
\omega_1+\omega_2 \quad & \mathbb{V}(z_1,z_2) \cup \mathbb{V}(z_1,z_4) \cup \mathbb{V}(z_4,2z_1z_3-z_2^2)& \frac{z_1^4z_4^2}{(2z_1z_3-z_2^2)^3}\\
\omega_1+\omega_3 \quad &\mathbb{V}(z_1,z_2) \cup \mathbb{V}(z_1,z_3) \cup \mathbb{V}(z_3,z_4)& \frac{(z_1z_4-z_2z_3)^2}{z_1z_3^3}\\
\omega_2+\omega_3 \quad &\mathbb{V}(z_1,z_2) \cup \mathbb{V}(z_2,z_3) \cup \mathbb{V}(z_1,2z_2z_4-z_3^2)& \frac{z_1^2(2z_2z_4-z_3^2)}{z_2^4}
\end{align*} 

It is easy to see that:

\begin{align*}
\omega_1+\omega_2&=z_1z_4(2z_1z_3-z_2^2) \Big(4\frac{dz_1}{z_1}+2\frac{dz_4}{z_4}-3\frac{d(2z_1z_3-z_2^2)}{2z_1z_3-z_2^2}\Big), \\
\omega_1+\omega_3&=z_1z_3(z_1z_4-z_2z_3)\Big(2\frac{d(z_1z_4-z_2z_3)}{z_1z_4-z_2z_3}-\frac{dz_1}{z_1}-3\frac{dz_3}{z_3}\Big),\\
\omega_2+\omega_3&=z_2z_1(2z_2z_4-z_3^2)\Big(4\frac{dz_2}{z_2}-2\frac{dz_1}{z_1}-\frac{d(2z_2z_4-z_3^2)}{2z_2z_4-z_3^2}\Big),\\
\end{align*} 

\noindent this means that these foliations are in the logarithmic component $\overline{L(1,1,2)}$. With this we obtain the following corollaries.

 \begin{cor} The exceptional component $\overline{E(3)}$ does not have foliations that have a line as singular set.
 \end{cor}
 
 \begin{cor} The boundary of the exceptional component is contained in a logarithmic component, more precisely:
 $$\overline{E(3)}-E(3) \subset L(1,1,2).$$
 \end{cor}

As we say before, the 1-forms $\omega_1$, $\omega_2$ and $\omega_3$ have a hyperplane as singular set, then they do not define a foliation but if we remove the singular hyperplane we obtain foliations which are linear pull-back of foliations on 
$\CP$ of degree $1$ and such that:

 \begin{align*}
 \textrm{1-form} \quad& \textrm{Singular Set} & \textrm{Rational First Integral}\\
\frac{\omega_1}{z_1} \quad & \mathbb{V}(z_1,z_3) \cup \mathbb{V}(z_1,z_4) \cup \mathbb{V}(z_3,z_4) &\frac{z_1z_4^2}{z_3^3}\\
\frac{\omega_2}{z_2} \quad &  \mathbb{V}(z_1,z_2) \cup \mathbb{V}(z_1,z_4) \cup \mathbb{V}(z_2,z_4)  &\frac{z_1^2z_4}{z_2^3}\\
\frac{\omega_3}{z_3} \quad &  \mathbb{V}(z_1,z_2) \cup \mathbb{V}(z_1,z_3) \cup \mathbb{V}(z_2,z_3) &\frac{z_2^2}{z_1z_3}\\
\end{align*}

\section{The logarithmic components of  $\mathcal{F}(2,3)$}

We know that there exist two components of $\mathcal{F}(2,3)$ which generic elements are logarithmic foliations. We are going to describe the singular set for the foliations in every component.
Remember that:

\begin{itemize}
\item  $\overline{L(1,1,1,1)}$ is the closure of the space of logarithmic foliations given by the 1-forms:

$$\omega=\sum_{i=1}^4 \lambda_i L_1L_2L_3L_4 \frac{dL_i}{L_i},$$

\noindent where $L_1, L_2, L_3$ and $L_4$ in $\C[z_1,z_2,z_3,z_4]$ have degree $1$ and $\sum_{i=1}^4 \lambda_i=0$.

\item $\overline{L(1,1,2)}$ is the closure of the space of logarithmic foliations given by the 1-forms:

$$\omega=\sum_{i=1}^3 \lambda_i f_1f_2f_3 \frac{df_i}{f_i},$$

\noindent where $f_1, f_2$ are different hyperplanes, $f_3$ is an irreducible hypersurface of degree $2$ and $\lambda_1+\lambda_2+2\lambda_3=0$.
\end{itemize}
\vspace{.3cm}

If $\omega= \sum_{i=1}^4 \lambda_i L_1L_2L_3L_4 \frac{dL_i}{L_i} \in L(1,1,1,1)$, using the part (c) of the proposition 2.1 in page 96 of \cite{Cerveau-Mattei} we have that:

$$Sing(\omega)=\bigcup_{i \neq j} \mathbb{V}(L_i,L_j),$$ 

\noindent and this is not a line because
$L_1, L_2, L_3$ and $L_4$ are different hyperplanes. On the other hand, with theorem 1.1 page 91 of \cite{Cerveau-Mattei} adapted to this case, we can prove that the foliations in $\overline{L(1,1,1,1)}-L(1,1,1,1)$ have the form:

\begin{enumerate}
\item $\omega_1=L_1^2L_2L_3\big(\sum_{i=1}^3 \lambda_i \frac{dL_i}{L_i}+d(\frac{\alpha}{L_1})\big)$,   where $\alpha$ is homogeneous of degree $1$.
\item $\omega_2=L_1^3L_2\big( \lambda_1 \frac{dL_1}{L_1}+\lambda_2 \frac{dL_2}{L_2}+d(\frac{\alpha}{L_1^2})\big)$, where $\alpha$ is homogeneous of degree $2$.
\item $\omega_3=L_1^4\big( \lambda_1 \frac{dL_1}{L_1}+d(\frac{\alpha}{L_1^3})\big)=\lambda_1L_1^3dL_1+L_1d\alpha-3\alpha dL_1$,  where $\alpha$ is homogeneous of degree $3$.
\end{enumerate}

\noindent By part (d) of the proposition 2.1 in page 96 of \cite{Cerveau-Mattei}, we obtain the singular set for every case

\begin{enumerate}
\item $Sing(\omega_1)= \mathbb{V}(L_1,L_2) \cup \mathbb{V}(L_1,L_3) \cup \mathbb{V}(L_2,L_3) \cup \mathbb{V}(L_1,\alpha)$
\item $Sing(\omega_2)= \mathbb{V}(L_1,L_2)  \cup \mathbb{V}(L_1,\alpha)$
\item $Sing(\omega_1)=  \mathbb{V}(L_1,\alpha)$
\end{enumerate}

If $\omega= \sum_{i=1}^3 \lambda_i f_1f_2f_3 \frac{df_i}{f_i} \in L(1,1,2)$, then using part (c) of proposition 2.1 in \cite{Cerveau-Mattei} we have that:

$$Sing(\omega)=\bigcup_{i \neq j} \mathbb{V}(f_i,f_j) \cup Sing(df_3),$$ 

\noindent which is not a line. Finally, we just have to note that the foliations on the border of this component have the same form as the foliations on the border of the above logarithmic component. Then in the logarithmic components
$\overline{L(1,1,1,1)}$ and $\overline{L(1,1,2)}$ we do not have foliations with a line as singular set.

\section{Foliations on $\mathbb{CP}^3$ of codimension 1 that have a line as singular set}

In this section we present the classification of codimension 1 foliations on $\mathbb{CP}^3$ of degree $2$ that have a line as singular set. 
\\

As we say in the introduction, the classification of these type of foliations is the question proposed on page 57 of the book  \cite{Deserti-Cerveau}, in this sense the authors find the unique $\mathcal{L}$-foliation on $\mathbb{CP}^3$ of degree $2$ with a line as singular set (see theorem 4.5 of \cite{Deserti-Cerveau}). For this reason it is convenient to begin this section by saying something about 
$\mathcal{L}$-foliations.

\begin{defin} Let $\mathcal{F}$ be a foliation on $\mathbb{CP}^3$ of codimension $1$. We will say that $\mathcal{F}$ is a $\mathcal{L}$-foliation of codimension $1$ if there exists a Lie subalgebra $\mathfrak{g}$ of the Lie algebra  of the group
$Aut(\mathbb{CP}^3)$, such that for a generic point $z \in \mathbb{CP}^3$ we have the following property:

$$\textrm{$\mathfrak{g}(z)$ is the tangent space to the leaf of $\mathcal{F}$ at $z$}.$$

\noindent In particular, we have that $\dim \mathfrak{g}(z)=\dim \{X(z): X \in \mathfrak{g}\}=2$.
\end{defin}

For example, from theorem 4.5 of \cite{Deserti-Cerveau} we have that foliations $\omega$ and $\omega_1+\omega_3$ in the exceptional component $\overline{E(3)}$  (see theorem \ref{excepcional}) are $\mathcal{L}$-foliations. With the same result we have also the following:

\begin{cor} \label{L-foliation} The foliation on $\mathbb{CP}^3$ of degree $2$ given by the 1-form:

$$\omega=3(z_1z_3^2+z_2z_3z_4+z_2^3)dz_3-z_3d(z_1z_3^2+z_2z_3z_4+z_2^3),$$ 

\noindent is the unique $\mathcal{L}$-foliation, up to linear change of coordinates, that have a line as singular set.

\end{cor}

Now we present the theorem of classification of these type of foliations.

 \begin{teo} \label{Classification} The foliations on $\mathbb{CP}^3$ of codimension $1$ of degree $2$ that have a line as singular set are, up to change of coordinates: 

\begin{enumerate}

\item The linear pull-back in $\mathbb{CP}^3$ of the foliations on $\CP$ of degree $2$ given by the 1-forms:

\begin{align*}
\nu_2&=-z_3^3dz_1+z_3(z_2^2+z_1z_3)dz_2+(z_1z_3^2-z_2^3-z_1z_2z_3)dz_3\\
\nu_3&=z_2^2z_3dz_1-z_3(z_2^2+z_3^2)dz_2+(z_2^3+z_2^2z_3-z_1z_2z_3)dz_3\\
\nu_4&=(z_1z_2z_3+z_3^3-z_2^3)dz_1+z_2(z_3^2-z_1z_2)dz_2-(z_2^2z_3+z_1^2z_2+z_1z_3^2)dz_3.
\end{align*}

\item Or the foliation given by the 1-form:
 $$\omega_f=3fdz_3-z_3df,$$ 

\noindent where the polynomial $f(z_1,z_2,z_3,z_4)$ is one of the following:

\begin{align*}
&z_1z_3^2+z_3z_4^2+z_2^3,\\
& z_1z_3^2+z_2z_3z_4+z_2^3, \\
& z_1z_3^2+z_2^3.
\end{align*}
 
 \noindent In these cases the associated foliation has the rational first integral $\frac{f}{z_3^3}$, and only the third one is a linear pull-back of a foliation on $\mathbb{CP}^2$.

\end{enumerate}
\end{teo}

\begin{proof} For the proof we will study every component of the space of foliations on $\mathbb{CP}^3$ of codimension 1 of degree $2$. As we mentioned before and following the results of \cite{Cerveau-LinsNeto}, we have 6 components in 
$\mathcal{F}(2,3)$ that we have called in the introduction: 1,...,6. Components 4, 5 and 6 were analyzed in the preceding sections where we saw that they do not have foliations that have a line as singular set. 
\\

For the study of component $\overline{R(2,2)}$ we proceed as follows: if the polynomials $f$ and $g$ define quadric surfaces in $\mathbb{CP}^3$, with $\mathbb{V}(f)$ smooth, then the polynomial $g$ is not a double line, because 
the degree of the foliation is $2$. 
Let $Q=\mathbb{V}(f)$, by adjuntion formula we have that $K_Q= \mathcal{O}_Q(2) \otimes K_{\mathbb{CP}^3}= \mathcal{O}_Q(-2)$, then $-K_Q=\mathcal{O}_Q(2)$, and this is the class in $Pic(Q)$ of the intersection with another reduced quadric.
\\

On the other hand, we know that $Q$ is isomorphic to $\mathbb{P}^1 \times \mathbb{P}^1$ then $K_Q= p_1^*K_{\mathbb{P}^1} \otimes p_2^*K_{\mathbb{P}^1}$, we conclude that $K_Q$ has class $(-2,-2)$.
Therefore $\mathcal{O}_Q(2)$ has class $(2,2) \neq (0,4)$, this says that the intersection of $Q$ with another reduced quadric surface is not a line. In the frontier of this component we have foliations with a rational first integral where one of the quadric is the product of two different hyperplanes, then we do not have a line as singular set for the associated foliation.
\\

Then all the foliations of $\mathcal{F}(2,3)$ that have a line as singular set are in the components 1 and 3. In the component 1 we have the linear pull-back of the foliations on $\CP$ of degree $2$ with a unique singularity. By \cite{Cerveau}
we know that the three that do not have rational first integral are $\nu_2$, $\nu_3$ and $\nu_4$.
\\

It remains study the component $\overline{R(1,3)}$: Let $\mathcal{F} \in \overline{R(1,3)}$ such that it has a line as singular set. We can suppose that the rational first integral is:

$$\frac{f(z_1,z_2,z_3,z_4)}{z_3^3}$$

\noindent where $f(z_1,z_2,z_3,z_4)$ defines a reduced cubic hypersurface in $\mathbb{CP}^3$. The algebraic variety $\mathbb{V}(f,z_3)$ is contained in the singular set $Sing(\mathcal{F})$ of the foliation, then we can suppose that this is $\mathbb{V}(z_2,z_3)$,
therefore 

$$f(z_1,z_2,z_3,z_4)=z_3h(z_1,z_2,z_3)+z_3z_4L(z_1,z_2,z_3,z_4)+az_2^3,$$ 

\noindent where $h(z_1,z_2,z_3)$ has degree 2, $L=a_1z_1+a_2z_2+a_3z_3+a_4z_4$, for some complex numbers $a_1, a_2, a_3, a_4$, and $a \in \C^*$.  
\\

Note that if $L=0$, then the foliation depends just of the variables $z_1,z_2,z_3$ and this is the linear pull-back 
of a foliation on $\CP$ of degree $2$ with a rational first integral and a unique singularity at $(1:0:0)$. 
\\

Up to change of coordinates there is only one foliation with this properties, and it has rational first integral:

$$\frac{z_1z_3^2+z_2^3}{z_3^3},$$

\noindent then we can take $h=z_1z_3$ and $a=1$. We can see that if $a_1 \neq 0$,  the cubic hypersurface defined by $z_1z_3^2+z_3z_4(\sum_{i=1}^4 a_iz_i)+z_2^3$ has a singularity out of the line $\mathbb{V}(z_2,z_3)$, 
since this is also a singularity for the foliation then we must ask that $a_1=0$. Then to finish the analysis we consider the following cases for $a_2, a_3$ and $a_4$:

\begin{enumerate}

\item Suppose that $a_4 \neq 0$, and consider $f=z_1z_3^2+z_3z_4(\sum_{i=2}^4 a_iz_i)+z_2^3,$ we can see that the unique singular point for the cubic hypersurface defined by $f$ is $(1:0:0:0)$. We have that, up to change of coordinates (see case
XX of the study of singular cubic hypersurfaces in section 9.2.3 of \cite{Dolgachev}), this hypersurface is defined by:

$$z_1z_3^2+z_3z_4^2+z_2^3.$$

\noindent This cubic hypersurface just contains the line $\mathbb{V}(z_2,z_3),$ which is the singular set of the foliation, in fact this is the unique cubic hypersurface which contains just one line (see table 9.1 of \cite{Dolgachev}). 
The foliation is not a linear pull-back of a foliation on $\mathbb{CP}^2$ because the hypersurface $\mathbb{V}(z_1z_3^2+z_3z_4^2+z_2^3)$ is not a cone over a cubic plane curve.
Finally, by theorem 4.5 of \cite{Deserti-Cerveau} we can conclude that this is not a $\mathcal{L}$-foliation either.

\item If $a_2 \neq 0$ and $a_4=0$, then the cubic hypersurface defined by 

$$f=z_1z_3^2+z_3z_4(a_2z_2+a_3z_3)+z_2^3,$$ 

\noindent has as singular set the line $\mathbb{V}(z_2,z_3)$ which is the singular set of the foliation. Then by the form of $f$ and theorem 9.2.1 of \cite{Dolgachev} this is, up to linear change of coordinates:

$$f=z_1z_3^2+z_2z_3z_4+z_2^3.$$ 

\noindent Hence the rational first integral of the foliation is $\frac{z_1z_3^2+z_2z_3z_4+z_2^3}{z_3^3}$ and this this the $\mathcal{L}$-foliation given in Corollary \ref{L-foliation}. Since this is a  $\mathcal{L}$-foliation and using proposition 3.7 in \cite{Deserti-Cerveau}  we have this is not a linear pull-back of a foliation on $\mathbb{CP}^2$. 

\item By parts 1 and 2, if $a_2 \neq 0$ or $a_4 \neq 0$ the foliation given by the 1-form $\omega_f$ is not a linear pull-back. If $a_2=a_4=0$ then the rational first integral of the foliation is:

$$\frac{(z_1+z_4)z_3^2+z_2^3}{z_3^3},$$

\noindent and with a linear change of coordinates we obtain that the rational first integral is equivalent to

$$\frac{z_1z_3^2+z_2^3}{z_3^3},$$

\noindent and the foliation is given by $\nu_1=z_3^3dz_1-z_1z_2^2dz_2+(z_2^3-z_1z_3^2)dz_2$, this is the linear pull-back of the unique foliation on $\mathbb{CP}^2$ of degree $2$ with a unique singular point and with a rational first integral. This means that this foliation is the unique in the intersection of the components $\overline{R(1,3)}$ and $S(2,3)$ with a line as singular set.
\end{enumerate}

\end{proof}

From the proof of the previous theorem we can conclude the following.

\begin{cor}  There exists, up to linear change of coordinates, only one codimension 1 foliation on $\mathbb{CP}^3$ of degree 2 with a line as singular set which is neither  $\mathcal{L}$-foliation nor  linear pull-back of a foliation
on $\mathbb{CP}^2$. 
\\

This foliation has the rational first integral $\frac{z_1z_3^2+z_3z_4^2+z_2^3}{z_3^3}$.
\end{cor}

We finish with the following examples of foliations in $\mathcal{F}(s,3)$, where $s \geq 3$, that have a line as singular set.

\begin{ex}  In order  to give a family of foliations on $\mathbb{CP}^3$ of codimension 1 of arbitrary degree $s$ that have a line as singular set and a rational first integral of degree $s+1$ we do a generalization of the construction that we did in the above theorem.
\\

Let $s \in \N$, $P(z_2,z_3)=\sum_{i=0}^s a_i z_2^i z_3^{s-i} \in \C[z_2,z_3]$ with $a_s \neq 0$, and $a \in \C$. We consider:

$$\frac{z_1z_3^{s}+a z_3z_4^s+Q(z_2,z_3)}{z_3^{s+1}},$$

 \noindent where $Q(z_2,z_3)=\sum_{i=0}^s \frac{a_i}{i+1}z_2^{i+1}z_3^{s-i}$. Then the 1-form :
 
 $$\omega_{a}= z_3^{s+1}dz_1 -z_3P(z_2,z_3) dz_2+(-z_1z_3^s+z_2P(z_2,z_3)+sa  z_3z_4^s)dz_3-s a z_3^2 z_4^{s-1} dz_4,$$
 
 \noindent defines a foliation on $\mathbb{CP}^3$ of degree $s$ and its singular set is the line $\mathbb{V}(z_2,z_3)$. This foliation is in the rational component $\overline{R(1,s+1)}$ of $\mathcal{F}(s,3)$. If $a=0$ then $\omega_0$ is the pull-back of the foliation: 
 
 $$\omega_{0,\CP}=z_3^{s+1}dz_1 -z_3P(z_2,z_3) dz_2+(-z_1z_3^s+z_2P(z_2,z_3))dz_3$$ 
 
 \noindent on $\CP$ with a unique singularity and with the rational first integral:

$$\frac{z_1z_3^s+Q(z_2,z_3)}{ z_3^{s+1}}.$$
\end{ex}

With the following example we show that for degree greater than 2 there exist logarithmic foliations with a line as singular set. For some cases they have a minimal rational first integral of degree not bounded. 

\begin{ex} Let $s_1, s_2, s_3 \in \N$, $a \in \C$ and 

\begin{align*}
f_a(z_1,z_2,z_3,z_4)&=z_2^{s_2}z_3^{s_3}(z_1z_2^{s_1-1}+az_4^{s_1})+z_2^{s_1+s_2+s_3}+z_3^{s_1+s_2+s_3}.\\
\end{align*}

\noindent Take $\lambda_1, \lambda_2, \lambda_3 \in \C^*$ such that $\lambda_1(s_1+s_2+s_3)+\lambda_2+\lambda_3=0$, then the 1-form

\begin{align*}
\omega_a=f_az_2z_3 \big(&\lambda_1 \frac{df_a}{f_a}+ \lambda_2 \frac{dz_2}{z_2}+ \lambda_3 \frac{dz_3}{z_3}\Big),
\end{align*}

\noindent defines a logarithmic foliation on $\mathbb{CP}^3$ of degree $s_1+s_2+s_3$, its singular set is the line $\mathbb{V}(z_2,z_3)$. These foliations are in the irreducible logarithmic component $\overline{L(1,1,s_1+s_2+s_3)}$ of $\mathcal{F}(s_1+s_2+s_3,3)$. If $ \lambda_2, \lambda_3 \in \N$, then 
the foliation has the minimal rational first integral:

$$\frac{f_a^{-\lambda_1}}{z_2^{\lambda_2} z_3^{\lambda_3}},$$

\noindent which has degree $\lambda_1+\lambda_2$, that means that the degree of this minimal rational first integral is not bounded. If $a=0$ the foliation is a linear pull-back of a foliation on $\mathbb{CP}^2$ of degree $s_1+s_2+s_3$.
\end{ex}

\end{document}